\documentclass[11pt]{amsart}
\usepackage{graphicx}
\usepackage{amssymb, amsmath}
\newtheorem{theorem}{Theorem}[section]
\newtheorem{corollary}[theorem]{Corollary}
\newtheorem{lemma}[theorem]{Lemma}

\theoremstyle{definition}

\theoremstyle{remark}

\begin{document}
\title{Polyadic groups and automorphisms of cyclic extensions}
\author{\sc M. Shahryari}
\thanks{{\scriptsize
\hskip -0.4 true cm MSC(2010): 20N15
}}

\address{ Department of Pure Mathematics,  Faculty of Mathematical
Sciences, University of Tabriz, Tabriz, Iran }
\email{mshahryari@tabrizu.ac.ir}
\date{\today}

%%% ----------------------------------------------------------------------
\begin{abstract}
We show that for any $n$-ary group $(G,f)$, the group $Aut(G,f)$ can be embedded in $Aut(\mathbb{Z}_{n-1}\ltimes G)$ and so we can obtain a class of interesting automorphisms of cyclic extensions.
\end{abstract}

\maketitle

%%%%%%%%%%%%%%%%%%%%%%%%%%%%%%%%%%%%%%%%%%%%%%%%%%%%%%%%%%%%%%%%%%%%%%
\section{Introduction}
In this article, we show that the automorphism group of any poyadic group  can be embedded in the group of automorphisms of its Post cover and then we apply this embedding  to obtain some interesting automorphisms of cyclic extensions. Our notations in this article are standard and can be find in \cite{Shah}, for example.

Let $(G,f)$ be an $n$-ary group. We know that, there is a binary operation (dot: $\cdot$) on $G$, such that $(G, \cdot)$ is an ordinary group, and further, there is a $\theta \in Aut(G, \cdot)$ with an element $b\in G$, such that\\
i- $\theta(b)=b$, and $\theta^{n-1}(x)=bxb^{-1}$, for all $x\in G$.\\
ii- $f(x_1^n)=x_1\theta(x_2)\cdots \theta^{n-1}(x_n)b$.\\
So, some times, we denote $(G, f)$  by the notation $der_{\theta, b}(G, \cdot)$. If $b=e$, the identity element of $(G, \cdot)$, then we use the notation $der_{\theta}(G, \cdot)$.\\

We associate another binary group to $(G,f)$ which is called the {\em universal covering group} or {\em Post's cover} of $(G,f)$. Let $a$ be an arbitrary element of $G$ and suppose $G^{\ast}_a=\mathbb{Z}_{n-1}\times G$. Define a binary operation on this set by
$$
(i,x)\ast (j,y)=(i+j+1, f_{\ast}(x, \stackrel{(i)}{a}, y, \stackrel{(j)}{a}, \overline{a}, \stackrel{(i, j)}{a})).
$$
Here of course, $i+j+1$ is computed modulo $n-1$, and $(i, j)=n-i-j-3,\ (mod\ n-1)$. The symbol $f_{\ast}$ indicates that $f$ applies one or two times depending on the values of $i$ and $j$ and $\overline{a}$ denotes the skew element of $a$. It is proved that (see \cite{Mich}), $G^{\ast}_a$ is binary group and the subset
$$
R=\{ (n-2, x): x\in G\}
$$
is a normal subgroup such that $G^{\ast}_a/R\cong \mathbb{Z}_{n-1}$. Further, if we identify $G$ by the subset
$$
\{ (0, x): x\in G\},
$$
then $G$ is a coset of $R$ and it generates $G^{\ast}_a$. We also have
$$
f(x_1^n)=x_1\ast x_2\ast \cdots \ast x_n.
$$
It is not hard to see that for all $a, b\in G$, we have $G^{\ast}_a\cong G^{\ast}_b$, so for simplicity, we always, assume that $a=e$, the identity element of $(G, \cdot)$.

Through this article, we  assume that $(G, f)=der_{\theta}(G, \cdot)$. So, we have $\theta^{n-1}=id$ and
$$
f(x_1^n)=x_1\theta(x_2)\cdots \theta^{n-2}(x_{n-1})x_n.
$$
We also assume that $e$ is the identity element of $(G, \cdot)$. We will prove first the following theorem on the structure of the Post's cover.

\begin{theorem}
We have $(der_{\theta}(G, \cdot))^{\ast}_e\cong \mathbb{Z}_{n-1}\ltimes G$, where $\mathbb{Z}_{n-1}$ acts on $(G, \cdot)$ by $i.x=\theta^i(x)$.
\end{theorem}
Note that, we used a special case of this theorem in \cite{Shah2}, to investigate representations of polyadic groups. The main idea of this article is almost the same as in \cite{Shah2}. Our second goal is to obtain an embedding from $Aut(G, f)$ to $Aut(G^{\ast}_e)$. The method we employ is the same as in \cite{Shah2}. For any $i\in \mathbb{Z}_{n-1}$ and $u\in G$, suppose $\delta(i,u)=\theta(u)\theta^2(u)\ldots \theta^i(u)$. We prove

\begin{theorem}
Let $\Lambda\in Aut(G, f)$ and define $\Lambda^{\ast}:G^{\ast}_e\to G^{\ast}_e$ by
$$
\Lambda^{\ast}(i, x)=(i, \Lambda(x)\delta(i, u)),
$$
where $u=\Lambda(e)$. Then the map $\Lambda\mapsto \Lambda^{\ast}$ is an embedding.
\end{theorem}

In \cite{Ham}, the structure of automorphisms of $(G, f)$ is determined. If $\Lambda\in Aut(G,f)$, then we have $\Lambda=R_u\phi$, where $u$ is an idempotent element, i.e. $f(\stackrel{(n)}{u})=u$,  $R_u$ is the right translation by $u$ and $\phi$ is an ordinary automorphism of $(G, \cdot)$ with the property $[\phi, \theta]=I_u$, (the bracket denotes the commutator $\phi\theta\phi^{-1}\theta^{-1}$ and $I_u$ is the inner automorphism corresponding to $u$). The converse is also true; if $u$ and $\phi$ satisfy above conditions, the $\Lambda=R_u\phi$ is an automorphism of the polyadic group $(G,f)$. We will use this fact frequently through this article. The interested reader should see \cite{Ham} for a full description of homomorphisms between polyadic groups.

Combining Theorems 1.1 and 1.2, we obtain an embedding of $Aut(G,f)$ into $Aut(\mathbb{Z}_{n-1}\ltimes G)$. More precisely, we prove the following.

\begin{theorem}
Let $\hat{G}=A\ltimes G$, with $A=\langle a\rangle$ cyclic of order $n-1$ and let $\theta(x)=axa^{-1}$. Then for any $\phi\in Aut(G)$ and $u\in G$, the hypotheses $[\phi, \theta]=I_u$ and $(au)^{n-1}=1$ imply that the map
$$
(a^i, x)\mapsto (a^i, u^{-1}\phi(x)u(au)^ia^{-i})
$$
is an automorphism of $\hat{G}$ and these automorphisms are mutually distinct.
\end{theorem}

\section{Proofs}
To prove 1.1, note that in $G^{\ast}_e$, we have
\begin{eqnarray*}
(i,x)\ast(j,y)&=&(i+j+1, f_{\ast}(x, \stackrel{(i)}{e}, y, \stackrel{(j)}{e}, \overline{e}, \stackrel{(i, j)}{e}))\\
              &=&(i+j+1, x\theta(e)\cdots\theta^i(e)\theta^{i+1}(y)\theta^{i+2}(e)\\
              &\ &\ \ \ \ \ \ \ \ \ \ \ \ \ \ \ \ \ \ \ \  \cdots\theta^{i+j+2}(e)\theta^{i+j+3}(\overline{e})\cdots \theta^{n-2}(e))\\
              &=&(i+j+1, x\theta^{i+1}(y)\theta^{i+j+3}(\overline{e})),
\end{eqnarray*}
but, since $\overline{e}=e$, so
\begin{eqnarray*}
(i,x)\ast(j,y)&=&(i+j+1, x\theta^{i+1}(y))\\
              &=&(i,x)(1,e)(j,y),
\end{eqnarray*}
where the right hand side product is done in $\mathbb{Z}_{n-1}\ltimes G$. Note that in general, if $(A, \cdot)$ is a group and $a\in A$, then we can define a new binary operation on $A$ by $x\circ y=xay$ and together with this new operation, $A$ is a group too, and so we denote it by $A_a=(A, \circ)$. We have $A\cong A_a$ and the isomorphism is given by $\varphi(x)=a^{-1}x$. Now, by this notation, we have
$$
(der_{\theta}(G, \cdot))^{\ast}_e=G^{\ast}_e=(\mathbb{Z}_{n-1}\ltimes G)_{(1,e)},
$$
and hence
$$
(der_{\theta}(G, \cdot))^{\ast}_e\cong\mathbb{Z}_{n-1}\ltimes G.
$$
This completes the proof of 1.1.\\

Now, let $\Lambda\in Aut(G, f)$ and $u=\Lambda(e)$. Define $\Lambda^{\ast}_e:G^{\ast}_e\to G^{\ast}_u$ by $\Lambda^{\ast}_e(i,x)=(i,\Lambda(x))$.

\begin{lemma}
$\Lambda^{\ast}_e$ is an isomorphism.
\end{lemma}

To prove this lemma, note that
\begin{eqnarray*}
\Lambda^{\ast}_e((i,x)\ast(j,y))&=&\Lambda^{\ast}_e(i+j+1, x\theta^{i+1}(y))\\
                                &=&(i+j+1, \Lambda(x\theta^{i+1}(y))).
\end{eqnarray*}
On the other hand,
\begin{eqnarray*}
\Lambda^{\ast}_e(i,x)\ast \Lambda^{\ast}_e(j,y)&=&(i, \Lambda(x))\ast(j, \Lambda(y))\\
                                               &=&(i+j+1, f_{\ast}(\Lambda(x), \stackrel{(i)}{u}, \Lambda(y), \stackrel{(j)}{u}, \overline{u}, \stackrel{(i, j)}{u})).
\end{eqnarray*}
But $f(\overline{u}, \stackrel{(n-1)}{u})=u$, so $\Lambda(f(v , \stackrel{(n-1)}{e}))=\Lambda(e)$, where $\Lambda (v)=\overline{u}$. Therefore $f(v , \stackrel{(n-1)}{e})=e$ and so $v=e$ and hence $\overline{u}=\Lambda(e)=u$. Now, we have
\begin{eqnarray*}
\Lambda^{\ast}_e(i,x)\ast \Lambda^{\ast}_e(j,y)&=&(i+j+1, \Lambda(f_{\ast}(x, \stackrel{(i)}{e}, y, \stackrel{(j)}{e}, e, \stackrel{(i, j)}{e})))\\
&=&(i+j+1, \Lambda(x\theta^{i+1}(y)))\\
&=&\Lambda^{\ast}_e((i,x)\ast(j,y)).
\end{eqnarray*}
This shows that $\Lambda^{\ast}_e$ is an isomorphism. \\

An element $u\in G$ is said to be {\em idempotent} if $f(\stackrel{(n)}{u})=u$. For an arbitrary element $u\in G$, we remember that the {\em right translation map} $R_u$ is defined by $R_u(x)=xu$. In \cite{Ham}, it is proved that every element of $Aut(G,f)$ can be uniquely represented as $R_u\phi$ with $u$ an idempotent and $\phi\in Aut(G, \cdot)$ satisfies $[\phi, \theta]=I_u$, where $I_u$ is the inner automorphism of $G$, corresponding to $u$. The converse is also true and so we have a complete description of automorphisms of $Aut(G, f)$ in terms of automorphisms of $(G, \cdot)$ and idempotents. Now, for any idempotent $u$ and $i\in \mathbb{Z}_{n-1}$, define
$$
\delta(i, u)=\theta(u)\theta^2(u)\cdots \theta^i(u).
$$
Note that for the case $i=0$, we have $\delta(0, u)=\delta(n-1, u)=e$. If $\Lambda\in Aut(G,f)$ and $u=\Lambda(e)$, then we define a map $q_u:G^{\ast}_u\to G^{\ast}_e$ by $q_u(i, x)=(i, x\delta(i, u))$.

\begin{lemma}
The map $q_u$ is an isomorphism.
\end{lemma}

To prove the lemma, we first assume that $i, j\neq 0$. Note that in $G^{\ast}_u$, we have
\begin{eqnarray*}
(i,x)\ast (j,y)&=&(i+j+1, f_{\ast}(x, \stackrel{(i)}{u}, y, \stackrel{(j)}{u}, \overline{u}, \stackrel{(i, j)}{u}))\\
               &=&(i+j+1, \Lambda(f_{\ast}(\Lambda^{-1}(x), \stackrel{(i)}{e}, \Lambda^{-1}(y), \stackrel{(j)}{e}, e, \stackrel{(i, j)}{e})))\\
               &=&(i+j+1, \Lambda(\Lambda^{-1}(x)\theta^{i+1}(\Lambda^{-1}(y)))).
\end{eqnarray*}
Now, as we said before, $\Lambda=R_u\phi$ such that $\phi\in Aut(G)$ and $[\phi, \theta]=I_u$. Therefore
\begin{eqnarray*}
(i,x)\ast (j,y)&=&(i+j+1, R_u\phi(\phi^{-1}R^{-1}_u(x)\theta^{i+1}(\phi^{-1}R^{-1}_u(y))))\\
               &=&(i+j+1, R_u((xu^{-1})\phi\theta^{i+1}\phi(yu^{-1}))).
\end{eqnarray*}
Since $[\phi, \theta]=I_u$, so we have $\phi\theta^{i+1}\phi^{-1}=(I_u\theta)^{i+1}$. But
\begin{eqnarray*}
(I_u\theta)^{i+1}(z)&=&u\theta(u)\cdots\theta^i(u)\theta^{i+1}(z)\theta(u)^{-1}\cdots\theta(u)^{-1}u^{-1}\\
                    &=&u\delta(i, u)\theta^{i+1}(z)\delta(i, u)^{-1}u^{-1}.
\end{eqnarray*}
Hence, we have
\begin{eqnarray*}
(i,x)\ast (j,y)&=&(i+j+1, R_u(xu^{-1}u\delta(i, u)\theta^{i+1}(yu^{-1})\delta(i, u)^{-1}u^{-1}))\\
               &=&(i+j+1, x\delta(i,u)\theta^{i+1}(yu^{-1})\delta(i,u)^{-1}).
\end{eqnarray*}
Now, we are ready to show that $q_u$ is a homomorphism. First, note that
\begin{eqnarray*}
q_u((i,x)\ast(j,y))&=&q_u(i+j+1, x\delta(i,u)\theta^{i+1}(yu^{-1})\delta(i,u)^{-1})\\
                   &=&(i+j+1, x\delta(i,u)\theta^{i+1}(yu^{-1})\delta(i,u)^{-1}\delta(i+j+1, u)).
\end{eqnarray*}
On the other hand
\begin{eqnarray*}
q_u(i,x)\ast q_u(j,y)&=&(i, x\delta(i,u))\ast(j, y\delta(j,u))\\
                     &=&(i+j+1, x\delta(i,u)\theta^{i+1}(y\delta(j,u)))\\
                     &=&(i+j+1, x\delta(i,u)\theta^{i+1}(y)\theta^{i+1}(\delta(j,u))).
\end{eqnarray*}
Hence, $q_u$ is a homomorphism, if and only if we have
$$
\theta^{i+1}(\delta(j,u))=\theta^{i+1}(u^{-1})\delta(i,u)^{-1}\delta(i+j+1, u).
$$
But, we have,
\begin{eqnarray*}
\theta^{i+1}(u^{-1})\delta(i,u)^{-1}\delta(i+j+1, u)&=&\theta^{i+2}(u)\cdots\theta^{i+j+1}(u)\\
                                                    &=&\theta^{i+1}(\delta(j,u)).
\end{eqnarray*}
The case $i=0$ can be verified similarly, so $q_u$ is a homomorphism. It is easy to see that also $q_u$ is a bijection and so we proved the lemma.\\

Combining two isomorphisms $q_u$ and $\Lambda^{\ast}_e$, we obtain an automorphism $\Lambda^{\ast}=q_u\circ \Lambda^{\ast}_e\in Aut(G^{\ast}_e)$. Note that, we have
$$
\Lambda^{\ast}(i,x)=(i, \Lambda(x)\delta(i,u))=(0,\Lambda(x))\ast(0,u)^i.
$$

\begin{lemma}
The map $\Lambda\mapsto \Lambda^{\ast}$ is an embedding from $Aut(G,f)$ into $Aut(G^{\ast}_e)$.
\end{lemma}

To verify this lemma, let $\Lambda_1, \Lambda_2\in Aut(G,f)$ and $u=\Lambda_1(e)$ and $v=\Lambda_2(e)$. Suppose also $w=\Lambda_1(v)=(\Lambda_1\circ\Lambda_2)(e)$. We have
$$
(\Lambda_1\circ\Lambda_2)^{\ast}(i,x)=(i, \Lambda_1(\Lambda_2(x))\delta(i, w)).
$$
On the other hand
\begin{eqnarray*}
\Lambda_1^{\ast}(\Lambda_2^{\ast}(i,x))&=&\Lambda_1^{\ast}(i, \Lambda_2(x)\delta(i,v))\\
                                     &=&(i, \Lambda_1(\Lambda_2(x)\delta(i, v))\delta(i,u)).
\end{eqnarray*}
But we have
\begin{eqnarray*}
\Lambda_1(\Lambda_2(x)\delta(i,v))&=&\Lambda_1(\Lambda_2(x)\theta(v)\cdots\theta^i(v)\theta^{i+1}(e)\cdots
                                        \theta^{n-2}(e)e)\\
                                  &=&\Lambda_1(f(\Lambda_2(x),\stackrel{(i)}{v}, \stackrel{(n-i-2)}{e},e))\\
                                  &=&f(\Lambda_1(\Lambda_2(x)), \stackrel{(i)}{w}, \stackrel{(n-i-2)}{u},\Lambda_1(e))\\
                                  &=&\Lambda_1(\Lambda_2(x))\delta(i,w)\theta^{i+1}(u)\cdots
                                  \theta^{n-2}(u)\Lambda_1(e).
\end{eqnarray*}
Note that we have
\begin{eqnarray*}
\theta^{i+1}\cdots\theta^{n-2}\Lambda_1(e)\delta(i,u)
                            &=&\theta^{i+1}(u)\cdots\theta^{n-2}(u)u\theta(u)\cdots\theta^i(u)\\
                            &=&e,
\end{eqnarray*}
because,
\begin{eqnarray*}
\theta(u)\cdots\theta^i(u)\theta^{i+1}(u)&\cdots&\theta^{n-2}(u)\Lambda_1(e)\\
    &=&u^{-1}\Lambda_1(f(\stackrel{(n)}{e}))\\
    &=&u^{-1}\Lambda_1(e)\\
    &=&e.
\end{eqnarray*}
Therefore we obtain
$$
\Lambda_1^{\ast}(\Lambda_2^{\ast}(i,x))=(i, \Lambda_1(\Lambda_2(x))\delta(i, w)),
$$
and this shows that the map $\Lambda\mapsto \Lambda^{\ast}$ is a homomorphism. Now suppose $\Lambda^{\ast}=id$. Then
$\Lambda(x)\delta(i, u)=x$ for all $x$ and $i$, so if we put $x=e$, then $\delta(i,u)=u^{-1}$ for all $i$. Assuming $i=1$, we get $\theta(u)=u^{-1}$ and so assuming $i=2$, we obtain $u^{-1}u=u^{-1}$, hence $u=e$ and consequently $\Lambda=id$. This completes the proof of the lemma.\\

Remember that we proved
$$
G^{\ast}_e=(\mathbb{Z}_{n-1}\ltimes G)_{(1,e)}\cong \mathbb{Z}_{n-1}\ltimes G,
$$
and this isomorphism is given by $\varphi(i,x)=(1,e)^{-1}(i,x)$. So,
\begin{eqnarray*}
\varphi(i,x)&=&(n-2, e)(i,x)\\
            &=&(n+i-2, \theta^{n-2}(x))\\
            &=&(i-1, \theta^{n-2}(x)).
\end{eqnarray*}
Now, for any $\Lambda\in Aut(G,f)$, define
$$
\alpha(\Lambda)=\varphi^{-1}\circ \Lambda^{\ast}\circ \varphi.
$$
Therefore $\alpha(\Lambda)$ is an automorphism of $\mathbb{Z}_{n-1}\ltimes G$ and the map $\Lambda\mapsto \alpha(\Lambda)$ is an embedding. We have
\begin{eqnarray*}
\alpha(\Lambda)(i,x)&=&\varphi^{-1}(i-1, \Lambda(\theta^{-1}(x))\delta(i-1, u))\\
                    &=&(1,e)(i-1, \Lambda(\theta^{-1}(x))\delta(i-1, u))\\
                    &=&(i, (\theta\Lambda\theta^{-1})(x)\theta(u^{-1})\delta(i,u)).
\end{eqnarray*}
Since $\Lambda=R_u\phi$, so $(\theta\Lambda\theta^{-1})(x)=(\theta\phi\theta^{-1})(x)\theta(u)$. Hence
$$
\alpha(\Lambda)(i,x)=(i, (\theta\phi\theta^{-1})(x)\delta(i,u)).
$$
On the other hand $\theta\phi\theta^{-1}=I_u^{-1}\phi$ and hence
$$
\alpha(\Lambda)(i,x)=(i, u^{-1}\phi(x)u\delta(i,u)).
$$
Summarizing, we obtain the following corollary;

\begin{corollary}
There is an embedding $\alpha:Aut(G,f)\to Aut(\mathbb{Z}_{n-1}\ltimes G)$, such that
$$
\alpha(\Lambda)(i,x)=(i, u^{-1}\phi(x)u\delta(i,u)).
$$
\end{corollary}

Now, we are ready to prove Theorem 1.3. Suppose $\hat{G}=A\ltimes G$ where $A=\langle a\rangle$ is a cyclic of order $n-1$. Define an automorphism of $G$ by $\theta(x)=axa^{-1}$, so $\theta^{n-1}=id$. Let
$$
(G,f)=der_{\theta}(G, \cdot).
$$
So, there is an embedding $\alpha:Aut(G,f)\to Aut(\hat{G})$ such that
$$
\alpha(\Lambda)(a^i,x)=(a^i, u^{-1}\phi(x)u\delta(i,u)).
$$
Since $u$ is an idempotent, so $f(\stackrel{(n)}{u})=u$, and therefore
$$
u\theta(u)\cdots\theta^{n-1}(u)=u,
$$
which implies that
$$
aua^{-1}a^2ua^{-2}\cdots a^{n-2}ua^{-(n-2)}u=e.
$$
Hence $(au)^{n-1}=1$. Similarly, $\delta(i,u)=(au)^ia^{-i}$, so for any $\phi\in Aut(G)$ and for any $u\in G$, the hypotheses
$$
(au)^{n-1}=1,\ \ [\phi, \theta]=I_u
$$
imply that the map
$$
(a^i, x)\mapsto (a^i, u^{-1}\phi(x)u(au)^ia^{-i})
$$
is an automorphism of $\hat{G}$. Clearly this is an embedding and hence the theorem is proved.

\end{document}